\documentclass[12pt]{article}

\usepackage{latexsym}
\usepackage{amsopn,amscd}
\usepackage{amsfonts}
\usepackage{amssymb}
\usepackage{amsthm}
\usepackage{amsmath}
\usepackage{a4wide}
\usepackage[square,authoryear]{natbib}

\newtheorem{Theorem}{Theorem}[section]
\newtheorem{Proposition}[Theorem]{Proposition}

\newtheorem{Example}[Theorem]{Example}

\def\empha{\em}

\def\pemphas{}

\def\Sigm{\mathcal S}

\numberwithin{equation}{section}

\begin{document}

\title{\bf The universal Hopf algebra
 associated with a Hopf-Lie-Rinehart algebra}

\author{
J.~Huebschmann
\\[0.3cm]
 USTL, UFR de Math\'ematiques\\
CNRS-UMR 8524
\\
59655 Villeneuve d'Ascq C\'edex, France\\
Johannes.Huebschmann@math.univ-lille1.fr
 }
\maketitle
\begin{abstract}
{We introduce a notion of Hopf-Lie-Rinehart algebra and show that
the universal algebra of a Hopf-Lie-Rinehart algebra acquires an
ordinary Hopf algebra structure.}
\end{abstract}

\noindent
\\[0.3cm]
\\[0.3cm]

\section{Introduction}

For a Lie algebra $\mathfrak g$, the universal enveloping algebra
$\mathrm U\mathfrak g$ is well known to acquire a cocommutative
{\em Hopf\/} algebra structure having $\mathfrak g$ as its module
of primitive generators. In this paper we will explore diagonal
maps for the universal algebra of a general {\em Lie-Rinehart\/}
algebra.

Within the framework of standard homological algebra, the diagonal
map on the universal enveloping algebra of an ordinary Lie algebra
induces the ring structure in Lie algebra cohomology. The ring
structure can be described entirely in terms of the {\em
Maurer-Cartan algebra\/} of the Lie algebra (the differential
graded algebra of alternating forms on the Lie algebra). The
terminology {\em Maurer-Cartan algebra\/} goes back at least to
\cite{vanesthr} and was prompted by the development of the subject
in the 1930's. The terminology
{\em Chevalley-Eilenberg algebra\/} is nowadays common as well.
As a side remark we note that the familiar Maurer-Cartan equation 
lives in a variant of the Maurer-Cartan algebra;
this equation, in turn,  has recently resurged
as an interesting topic in its own, in particular in the theory of 
deformations and as the master equation in physics.
It is also worthwhile recalling that in view of a
classical result, in characteristic zero, any cocommutative Hopf
algebra is the enveloping algebra of its Lie algebra of primitive
elements, and this fact generalizes even to the differential
graded setting.

The question under what circumstances
the universal algebra
$\mathrm U(A,L)$ of a general {\em Lie-Rinehart\/} algebra $(A,L)$
has a diagonal map 
turning this universal algebra into a Hopf algebra has hardly been
explored in the literature, though. When $(A,L)$ is the
Lie-Rinehart algebra of smooth functions $A=C^{\infty}(N)$ and
smooth vector fields $L=\mathrm{Vect}(N)$ on a smooth manifold
$N$, the algebra $\mathrm U(A,L)$ is the algebra of (globally defined)
differential operators on $N$. The de Rham cohomology of $N$ then
amounts to the appropriate Ext-functors over $\mathrm U(A,L)$
\cite{rinehone}. For a general smooth manifold
$N$, there is no obvious way to put a
diagonal map on $\mathrm U(A,L)$, and the ring structure in
cohomology is defined directly in terms of the corresponding {\em
Maurer-Cartan algebra\/} (the differential graded algebra of de
Rham forms on the manifold). As a side remark we note that this
Maurer-Cartan algebra is the starting point for higher homotopies
generalizations of the structure, worked out in \cite{quasi}. This
kind of generalization occurs in nature, e.~g. in the theory of
foliations.

When $A$ is an ordinary Hopf algebra and when $\mathfrak g$ is an
ordinary Lie algebra acting on $A$ by derivations compatibly with
the diagonal, the crossed product algebra $A\odot \mathrm U\mathfrak g$
inherits a Hopf algebra structure in an obvious manner. This
situation arises e.~g. when $A$ is  the {\em Hopf\/} algebra of
algebraic functions on an algebraic group $H$ (e.g. a compact Lie
group) and when $\mathfrak g$ is the Lie algebra of an algebraic
group $G$ acting on $H$ by group automorphisms. 
The resulting Hopf algebra is no longer
cocommutative, though, unless $H$ is abelian.

Prompted by a recent posting \cite{moermrcu}, we decided to
communicate some structural insight for the case of a general
Lie-Rinehart algebra which we have been familiar with for many
years. We have already noted that the universal enveloping algebra
of an ordinary Lie algebra acquires a Hopf algebra structure in an
obvious way. The point we wish to emphasize here is that in order
for the naive extension of the fact just quoted to be valid for a
general Lie-Rinehart algebra $(A,L)$, a Hopf-algebra structure on
$A$ is needed. This is certainly consistent with the case of an
ordinary Lie algebra, the ground ring or ground field being
endowed with its obvious Hopf algebra structure (where all
structure maps come down to the identity).

We will introduce what we will refer to as a 
{\em
bi-Lie-Rinehart algebra\/} and we will, furthermore, 
refine this notion to that of a
{\em
Hopf-Lie-Rinehart algebra\/}, the crossed product algebra relative
to a Lie group being a special case. 
We shall show that the universal algebra of a bi- and, likewise,
 that of a Hopf-Lie-Rinehart algebra,
acquires a comultiplication and a counit which turn the universal algebra into a
bialgebra  or  Hopf algebra as appropriate.

I am indebted to Jim Stasheff for a number of comments
which helped improve the exposition.

\section{Lie-Rinehart algebras}

Let $R$ be a commutative ring, fixed throughout; the unadorned
tensor product symbol $\otimes$ will always refer to the tensor
product over $R$. Further, let $A$ be a commutative $R$-algebra,
let $L$ be an $A$-module (the action being written as
${(a\otimes \alpha) \mapsto a\alpha}$)  which is also an $R$-Lie algebra,
and suppose that $L$ acts on the left of $A$ by derivations
(the
action being  written as ${(\alpha \otimes a) \mapsto \alpha a}$).
Following {\em Rinehart\/} \cite{rinehone}, 
we will refer to $L$ as an $(R,A)$-{\em Lie algebra\/} 
provided suitable compatibility conditions are satisfied which generalize
standard properties of the Lie algebra of vector fields on a
smooth manifold viewed as a module over its ring of functions;
these conditions read
\begin{align} (a\,\alpha)(b) &= a\,(\alpha(b)), \quad \alpha \in L,
\,a,b \in A, \label{1.1.a}
\\
\lbrack \alpha, a\,\beta \rbrack &= a\, \lbrack \alpha, \beta
\rbrack + \alpha(a)\,\beta,\quad \alpha,\beta \in L, \,a \in A.
\label{1.1.b}
\end{align}
Occasionally we will spell out the $L$-action on $A$ explicitly in
the form $\omega\colon L \to \mathrm{Der}(A)$, so that
$(\omega(\alpha))(a)=\alpha(a)$. When the emphasis is on the pair
$(A,L)$, with the mutual structure of interaction, we refer to 
$(A,L)$ as a
{\em Lie-Rinehart algebra\/}. Given two Lie-Rinehart algebras
$(A,L)$ and $(A',L')$, a {\em  morphism \/} ${ (\varphi,\psi)
\colon (A,L) \longrightarrow (A',L') }$ of {\em  Lie-Rinehart
algebras\/} is the obvious thing, that is $\varphi\colon A \to A'$ 
is a morphism of $R$-algebras,
$\psi\colon L \to L'$ is a morphism of $R$-Lie algebras,
and these morphisms are compatible
with the additional structure. 
More precisely, the obvious diagrams
\begin{equation*}
\begin{CD}
A\otimes L
@>>>
A'\otimes L' @. \phantom{aaaaaa} L\otimes A @>>> L'\otimes A'
\\
@VVV
@VVV
 \phantom{aaaaaa}
@VVV
@VVV
\\
L @>>> L' @. \phantom{aaaaaa} A @>>> A'
\end{CD}
\end{equation*}
are commutative, the unlabelled vertical arrows being the 
corresponding
structure maps.
With this notion of morphism,
Lie-Rinehart algebras constitute a category. Apart from the
example of smooth functions and smooth vector fields on a smooth
manifold, a related (but more general) example is the pair
consisting of a commutative algebra $A$ and the $R$-module
$\mathrm {Der}(A)$ of derivations of $A$ with the obvious
$A$-module structure; here the commutativity of $A$ is crucial.

Given an $(R,A)$-Lie algebra $L$, its {\it universal algebra\/}
${(\mathrm U(A,L),\iota_L,\iota_A)}$ is an $R$-algebra $\mathrm
U(A,L)$ together with a morphism ${ \iota_A \colon A
\longrightarrow \mathrm U(A,L) }$ of $R$-algebras and a morphism
${ \iota_L \colon L \longrightarrow \mathrm U(A,L) }$ of Lie
algebras over $R$ having the properties
\begin{equation*}
\iota_A(a)\iota_L(\alpha) = \iota_L(a\,\alpha),\quad
\iota_L(\alpha)\iota_A(a) - \iota_A(a)\iota_L(\alpha) =
\iota_A(\alpha(a)),
\end{equation*}
and ${(\mathrm U(A,L),\iota_L,\iota_A)}$ is {\em universal\/}
among triples ${(B,\varphi_L,\varphi_A)}$ having these properties.
For example, when
 $A$ is the algebra of smooth functions on a smooth
manifold $N$ and  $L$  the Lie algebra of smooth vector fields on
$N$, then $\mathrm U(A,L)$ is the {\em algebra of (globally
defined) differential operators on\/} $N$. An explicit
construction for the $R$-algebra ${\mathrm U(A,L)}$ is given in
\cite{rinehone}. See our paper \cite{poiscoho} for an alternate
construction which employs the {\em Massey-Peterson\/}
\cite{masspete}  algebra.

The universal algebra $\mathrm U(A,L)$ admits an obvious filtered
algebra structure 
\begin{equation}
\mathrm U_{-1} \subseteq \mathrm U_{0}\subseteq
\mathrm U_{1}\subseteq \dots, 
\label{filt}
\end{equation}
cf. \cite{rinehone}, where
$\mathrm U_{-1}(A,L) = 0$ and where, for $p \geq 0$, $\mathrm
U_p(A,L)$ is the left $A$-submodule of $\mathrm U(A,L)$ generated
by products of at most $p$ elements of the image ${\overline L}$
of $L$ in $\mathrm U(A,L)$, and the associated graded object
$E^0(\mathrm U(A,L))$ inherits a commutative graded $A$-algebra
structure. We will refer to the filtration
\eqref{filt} as the {\em Poincar\'e-Birkhoff-Witt\/} filtration.
The Poincar\'e-Birkhoff-Witt Theorem for $\mathrm
U(A,L)$ then takes the following form where $\Sigm_A[L]$ denotes
the symmetric $A$-algebra on $L$, cf. (3.1) of \cite{rinehone}.

\begin{Theorem}[Rinehart]
\label{1.1} For an $(R,A)$-Lie algebra $L$ which is projective as
an $A$-module,
 the canonical
$A$-epimorphism $\Sigm_A[L] \longrightarrow \mathrm E^0(\mathrm
U(A,L)) $ onto the associated graded $A$-algebra  $\mathrm
E^0(\mathrm U(A,L)) $ is an isomorphism of $A$-algebras.
\end{Theorem}

Consequently, for an $(R,A)$-Lie algebra $L$ which is projective
as an $A$-module,
 the morphism
$\iota_L \colon L \longrightarrow \mathrm U(A,L) $ is injective,
and we can refer to $\mathrm U(A,L) $ as an {\em enveloping
algebra\/}.

It is worthwhile noting that, for an ordinary Lie algebra 
$\mathfrak g$ over a field of characteristic zero, 
one way to prove the 
Poincar\'e-Birkhoff-Witt Theorem consists in noting that the
coalgebra structure which underlies the obvious Hopf algebra structure
of $\mathrm U\mathfrak g$ is that of the symmetric coalgebra cogenerated by
$\mathfrak g$. This kind of reasoning breaks down for a general Lie-Rinehart
algebra $(A,L)$ since there is in general no obvious way to endow
the universal algebra $\mathrm U(A,L)$ with a coalgebra structure having the
correct features. Indeed, in this paper, we will spell out the requisite
additional Lie-Rinehart structure to arrive at the desired coalgebra structure
on  $\mathrm U(A,L)$.

\section{Induced Lie-Rinehart structures}

The notion of induced Lie-Rinehart structure 
 been introduced in \cite{poiscoho}. A geometric
analogue thereof can be found in \cite{higgmack}. For
intelligibility, we will now recall these induced structures:

Let $A$ and $A'$ be commutative $R$-algebras, let $L$ be an
$(R,A)$-Lie algebra, with structure maps $\omega \colon L
\longrightarrow \mathrm{Der}(A) $ and $[\,\cdot \,  , \, \cdot \,]
\colon L \otimes _R L \to L$, let $\tilde\omega \colon L
\longrightarrow \mathrm{Der}(A') $ be an action of $L$, viewed as
an $R$-Lie algebra, on $A'$ by derivations (but $L$ is not assumed
to admit an $A'$-module structure), and let ${\varphi \colon A
\longrightarrow A' }$ be a morphism of algebras which is also a
morphism of $L$-modules. Under these circumstances, let $L' = A'
\otimes _A L$, the induced $A'$-module. Consider the obvious
pairings
\begin{equation}
A' \otimes_R L \otimes _R A' \otimes_R L \longrightarrow A'
\otimes _A L, \label{1.16.1}
\end{equation}
given by $u \otimes \alpha \otimes v \otimes \beta \longmapsto uv
\otimes \lbrack \alpha,\beta \rbrack -(v\beta(u))\otimes \alpha +
(u \alpha(v)) \otimes \beta,$ where $u,v \in A',\,\alpha,\beta \in
L,$ and
\begin{equation}
A' \otimes_R L \otimes_R A' \longrightarrow A', \quad u \otimes
\alpha \otimes v \longmapsto u\cdot \alpha(v), \quad u,v \in
A',\alpha \in L. \label{1.16.2}
\end{equation}

A straightforward verification establishes the following, cf.
Proposition 1.16 in \cite{poiscoho}.

\begin{Proposition}
\label{1.16} Suppose that the action $\widetilde \omega$ of $L$ on
$A'$ is a morphism of $A$-modules where $\mathrm{Der}(A')$ is
turned into an $A$-module via $\varphi$. Then {\rm \eqref
{1.16.1}} induces an $R$-Lie algebra structure
$
[\, \cdot \, , \, \cdot \,]' \colon L' \otimes _R L'
\longrightarrow L'
$
on $L'$, {\rm \eqref {1.16.2}} induces an action $\omega' \colon
L' \to \mathrm{Der}(A')$ of $L'$ on $A'$ by derivations, and $[\,
\cdot \, , \, \cdot \,]'$ and $\omega'$ endow $L'$ with an
${(R,A')}$-Lie algebra structure in such a way that
\[
(\varphi,\varphi \otimes \mathrm{Id}) \colon (A,L) \longrightarrow
(A',L')
\]
is a morphism of Lie-Rinehart algebras.
\end{Proposition}

In the situation of Proposition \ref{1.16} we refer to $(L',[\,
\cdot \, , \, \cdot \,]',\omega')$ as the $(R,A')$-Lie algebra
 {\em induced from\/}
${\varphi}$; often we will then write $L'$ rather than $(L',[\,
\cdot \, , \, \cdot \,]',\omega')$.

\begin{Example} \label{ex} Let $A$ be a commutative $R$-algebra, let
$\mathfrak g$ be an $R$-Lie algebra, and let $\omega \colon
\mathfrak g \to \mathrm{Der}(A)$ be an action of $\mathfrak g$ on
$A$ by derivations. Then with the obvious change in notation, {\rm
\eqref{1.16.1}} and {\rm \eqref{1.16.2}} endow $A \otimes
\mathfrak g$ with an $(R,A)$-Lie algebra structure, referred to as
the {\em crossed product\/} of $A$ and $\mathfrak g$, cf. {\rm
\cite{mallia}}, and we will use the notation $A \odot \mathfrak g$
for the crossed product $(R,A)$-Lie algebra. In particular, over
the reals $\mathbb R$ as ground ring, given a real algebraic Lie
group $G$, with Lie algebra $\mathfrak g$, let $A$ be the algebra
of algebraic functions on $G$, endowed with the obvious $\mathfrak
g$-action. Then the crossed product $(\mathbb R,A)$-Lie algebra $A
\odot \mathfrak g$ amounts to the Lie algebra of algebraic vector
fields on $G$. The smooth analogue of this construction yields the
ordinary Lie algebra of smooth vector fields on $G$.
\end{Example}

\section{Bi-Lie-Rinehart algebras}

Let $A$ be a commutative $R$-bialgebra, not necessarily cocommutative,
with comultiplication $\Delta\colon A \to A\otimes A$ and counit
$\varepsilon\colon A \to R$, and let $L$ be an
$(R,A)$-Lie algebra, with structure maps $\omega \colon L
\longrightarrow \mathrm{Der}(A) $ and $[\,\cdot \,  , \, \cdot \,]
\colon L \otimes _R L \to L$. Let $\omega^{\otimes}\colon L \to
\mathrm{Der}(A \otimes A)$ be an action of  $L$, viewed as an
$R$-Lie algebra, on the tensor product algebra $A \otimes A$ by
derivations. To simplify the notation, we will then write
\begin{equation}
\alpha(a\otimes b) = (\omega^{\otimes} (\alpha))(a\otimes b).
\label{ob2}
\end{equation}
Thus, letting $A'=A\otimes A$, $\varphi=\Delta$, and $\widetilde
\omega = \omega^{\otimes}$, we have the data needed to endow the
induced $(A\otimes A)$-module $L^{\otimes}=(A\otimes A) \otimes_A
L$ with an $(R,A\otimes A)$-Lie algebra structure. In view of
Proposition \ref{1.16}, the following is immediate.

\begin{Proposition}
Suppose that, for every $a,b,c \in A$ and for every $\alpha \in
L$,
\begin{equation}
(a\alpha)(b\otimes c) = (\Delta(a)) \alpha(b\otimes c).
\label{ob3}
\end{equation}
Then {\rm \eqref
{1.16.1}} induces  an $R$-Lie algebra structure $[\, \cdot \, , \,
\cdot \,]^{\otimes} \colon L^{\otimes} \otimes _R L^{\otimes} \to
L^{\otimes}$ on $L^{\otimes}$, {\rm \eqref {1.16.2}} induces an
action $\omega^{\otimes} \colon L^{\otimes} \to \mathrm{Der}(A
\otimes A)$ of $L^{\otimes}$ on the tensor product algebra
 $A \otimes A$ by derivations, and $[\, \cdot \, , \, \cdot \,]^{\otimes}$ and
$\omega^{\otimes}$ endow $L^{\otimes}$ with an ${(R,A \otimes
A)}$-Lie algebra structure in such a way that
\begin{equation}
(\Delta,\Delta \otimes \mathrm{Id}) \colon (A,L) \longrightarrow
(A\otimes A,L^{\otimes}) \label{ob4}
\end{equation}
is a morphism of Lie-Rinehart algebras.
\end{Proposition}

Under the circumstances spelled out just before the previous 
proposition, we will say that
$(A,\Delta,L,\omega^{\otimes})$ is  a {\em bi-Lie-Rinehart
algebra\/} provided $(A,\Delta,L,\omega^{\otimes})$ satisfies the
requirement \eqref{ob3}
and, furthermore,
\begin{equation*}
(\varepsilon,0) \colon (A,L) \longrightarrow (R,0)
\end{equation*}
is a morphism of Lie-Rinehart algebras. 
We use the terminology  {\em bi-Lie-Rinehart
algebra\/} to avoid conflict with the notion of
 {\em Lie-Rinehart
bialgebra\/} which, in turn, has been introduced in \cite{banach}
as an abstraction from the notion of Lie bialgebroid.
We shall come back to Lie-Rinehart bialgebras in the last section.

Given a bialgebra $A$ and an $A$-module $L$, the symmetric
$A$-algebra $\Sigm_A[L]$ acquires an $R$-bialgebra structure in
an obvious way, the submodule $L$ being primitive.

\begin{Theorem} \label{th1}
The universal algebra $\mathrm U(A,L)$ of a bi-Lie-Rinehart
algebra acquires a comultiplication
\begin{equation}
\Delta\colon \mathrm U(A,L) \longrightarrow \mathrm U(A,L)\otimes
\mathrm U(A,L) \label{diag}
\end{equation}
and counit
\begin{equation}
\varepsilon\colon \mathrm U(A,L) \longrightarrow R \label{coun}
\end{equation}
turning $\mathrm U(A,L)$ into an $R$-bialgebra. Furthermore, the
graded algebra  $\mathrm E^0(\mathrm U(A,L)) $ associated
with the Poincar\'e-Birkhoff-Witt filtration
acquires
an obvious $R$-bialgebra structure and, when $L$ is projective
as an $A$-module, the canonical $A$-epimorphism 
\[
\Sigm_A[L]
\longrightarrow \mathrm E^0(\mathrm U(A,L)) 
\]
 onto  $\mathrm
E^0(\mathrm U(A,L)) $ is an isomorphism of $R$-bialgebras.
\end{Theorem}

\begin{proof} The direct sum $L\otimes A \oplus A \otimes L$
acquires an obvious $(R,A\otimes A)$-Lie algebra structure and the
standard $R$-Lie algebra diagonal map $L \to L \oplus L$ 
which assigns the pair $(x,x)$ to $x\in L$
induces a morphism
\begin{equation}
L^{\otimes} \longrightarrow L\otimes A \oplus A \otimes L
\label{ob5}\end{equation} of $(R,A\otimes A)$-Lie algebras. Thus
the composite of \eqref{ob4} and \eqref{ob5} yields a morphism
\begin{equation}
(\Delta,\Delta) \colon (A,L) \longrightarrow (A\otimes A,L\otimes
A \oplus A \otimes L) \label{ob6}
\end{equation}
of Lie-Rinehart algebras. This morphism, in turn, induces the
morphism
\begin{equation}
\mathrm U(\Delta,\Delta) \colon \mathrm U(A,L) \longrightarrow
\mathrm U(A\otimes A,L\otimes A \oplus A \otimes L) \label{ob7}
\end{equation}
between the universal algebras. It remains to show that the
algebra $\mathrm U(A,L)\otimes\mathrm U(A,L)$ is canonically
isomorphic to the universal algebra $\mathrm U(A\otimes A,L\otimes
A \oplus A \otimes L)$. In order to justify this claim, it
suffices to note that $\mathrm U(A,L)\otimes\mathrm U(A,L)$
satisfies the corresponding universal property. Indeed, relative
to the obvious morphism
\[
 \iota_{A\otimes A}= \iota_{A}\otimes \iota_{A}\colon A\otimes A \longrightarrow \mathrm U(A,L)
 \otimes \mathrm U(A,L)
\]
of $R$-algebras and relative to the obvious morphism
\[
\iota_{L\otimes A \oplus A \otimes L} =\iota_{L}\otimes
\iota_{A}\oplus \iota_{A}\otimes \iota_{L}
\colon L\otimes A
\oplus A \otimes L \longrightarrow \mathrm U(A,L) \otimes \mathrm
U(A,L)
\]
of Lie algebras over $R$, the algebra $\mathrm
U(A,L)\otimes\mathrm U(A,L)$ plainly satisfies the universal
property which characterizes the universal algebra $\mathrm
U(A\otimes A,L\otimes A \oplus A \otimes L)$.

By assumption, $(\varepsilon,0) \colon (A,L) \longrightarrow (R,0)$
is a morphism of Lie-Rinehart algebras. The induced morphism
\begin{equation*}
\varepsilon\colon\mathrm U(A,L) \longrightarrow \mathrm U(R,0) =R
\end{equation*}
of algebras yields the requisite counit where the notation
$\varepsilon$ is abused somewhat.
This establishes the theorem.
\end{proof}

\section{Hopf-Lie-Rinehart algebras}

Let $(A,L)$ be a Lie-Rinehart algebra.
Let $L^-$ be the $R$-Lie algebra {\em opposite\/} to $L$, that is,
as an $R$-module, $L^-$ coincides with $L$ whereas the bracket
on  $L^-$ is the negative of the bracket on $L$.

The $L$-action on $A$ by derivations being written as
$\omega\colon L \to \mathrm{Der}(A)$, let
\begin{equation*}
\omega^- = -\omega \colon L^- \to \mathrm{Der}(A).
\end{equation*}
This is plainly an  $L^-$-action on $A$ by derivations and,
together with the $A$-module structure on $L^-$ (which, as an $R$-module,
 coincides with $L$), the $A$-module $L^-$ thus acquires
an $(R,A)$-Lie algebra structure. In other words,
$(A,L^-)$ is a Lie-Rinehart algebra, which we will refer to as
the  {\em Lie-Rinehart algebra opposite to\/} $(A,L)$.
We remind the reader that, given an $R$-algebra $\mathrm U$, 
the {\em opposite\/} algebra $\mathrm U^{\mathrm{opp}}$ 
has the same underlying $R$-module as $\mathrm U$, with multiplication 
given by $x^{\mathrm{opp}}y^{\mathrm{opp}}=(yx)^{\mathrm{opp}}$
($x,y \in \mathrm U$). The canonical $R$-Lie algebra morphism from $L$ to 
 $\mathrm U(A,L)$ is as well an $R$-Lie algebra morphism from $L^-$ to
$\mathrm U(A,L)^{\mathrm{opp}}$ and, indeed,
$\mathrm U(A,L)^{\mathrm{opp}}$ satisfies the corresponding universal property
so that
$\mathrm U(A,L)^{\mathrm{opp}}$,
together with the obvious morphisms
$L^{-} \to \mathrm U(A,L)^{\mathrm{opp}}$ 
(which, as a morphism of $A$-modules, is just the morphism
$L \to \mathrm U(A,L)$
associated with the universal algebra for $(A,L)$)
and
$A \to \mathrm U(A,L)^{\mathrm{opp}}$,
yields the universal algebra
$\mathrm U(A,L^{-})$.

Let $(\Delta,\varepsilon)$ be
a bialgebra structure on $A$, and suppose that 
$(A,\Delta,L,\omega^{\otimes})$ is  a {\em bi-Lie-Rinehart
algebra\/}. Furthermore, let $S\colon A \to S$ be an {\em antipode\/}
turning $(A,\Delta,\varepsilon)$ into a {\em Hopf\/}-algebra.
We will refer to
$(A,\Delta,S,L,\omega^{\otimes})$ as  a {\em Hopf-Lie-Rinehart
algebra\/} provided 
\begin{equation}
(S, -) \colon (A,L)\longrightarrow (A,L^-)
\label{anti3}
\end{equation}
is a morphism of Lie-Rinehart algebras;
here $- \colon L\to L^-$ refers to the $A$-linear map which sends
$x\in L$ to $-x$, the   $L^-$ underlying $A$-module
being the same as that underlying $L$.

\begin{Theorem} \label{th2}
The universal algebra $\mathrm U(A,L)$ of a Hopf-Lie-Rinehart
algebra acquires a comultiplication
\begin{equation}
\Delta\colon \mathrm U(A,L) \longrightarrow \mathrm U(A,L)\otimes
\mathrm U(A,L), \label{diag2}
\end{equation}
counit
\begin{equation}
\varepsilon\colon \mathrm U(A,L) \longrightarrow R, \label{coun2}
\end{equation}
and antipode
\begin{equation}
S\colon \mathrm U(A,L) \longrightarrow \mathrm U(A,L) \label{anti2}
\end{equation}
turning $\mathrm U(A,L)$ into an $R$-Hopf algebra. Furthermore, the
associated graded algebra  $\mathrm E^0(\mathrm U(A,L)) $ acquires
an obvious $R$-Hopf algebra structure and, when $L$ is projective
as an $A$-module, the canonical $A$-epimorphism $\Sigm_A[L]
\longrightarrow \mathrm E^0(\mathrm U(A,L)) $ onto  $\mathrm
E^0(\mathrm U(A,L)) $ is an isomorphism of $R$-Hopf algebras.
\end{Theorem}

\begin{proof}
In view of Theorem \ref{th1}, it remains to establish the existence of 
the antipode. However, this is straightforward:
By definition, the antipode $S$ fits into a morphism of Lie-Rinehart algebras
of the kind \eqref{anti3} and this morphism in turn, induces the morphism
\begin{equation*}
\mathrm U(S, -) \colon \mathrm U(A,L)\longrightarrow \mathrm U(A,L^-)
\end{equation*}
of $R$-algebras. However, as $R$-modules,
$ \mathrm U(A,L^-)$ and $\mathrm U(A,L)$ coincide whence 
$\mathrm U(S, -)$ yields the requisite antipode.
We leave the details to the reader.
\end{proof}

\begin{Example}
Let $(A,L)$ be a Lie-Rinehart algebra and suppose that, for some
$R$-Lie algebra $\mathfrak g$ acting on $A$ by derivations, the
$(R,A)$-Lie algebra $L$ can be written as the crossed product
$(R,A)$-Lie algebra $A\odot \mathfrak g$. Furthermore, suppose
that $\Delta \colon A \to A \otimes A$ and
$\varepsilon \colon A \to R$
constitute a coalgebra structure turning
$A$ into a Hopf algebra such that, relative to the obvious
$\mathfrak g$-action on $A\otimes A$, $\Delta$ and $\varepsilon$ are
 compatible with
the $\mathfrak g$-actions. Then the crossed product $(A\otimes A)
\odot \mathfrak g$ is an $(R,A\otimes A)$-Lie algebra, and the
comultiplication $\Delta$ and counit $\varepsilon$ induce  morphisms
\[
(A,A\odot \mathfrak g) \longrightarrow (A\otimes A ,(A\otimes A)
\odot \mathfrak g)
\]
and
\begin{equation*}
(\varepsilon,0) \colon (A,L) \longrightarrow (R,0)
\end{equation*}
of Lie-Rinehart algebras. This yields the requisite data to apply
Theorem {\rm \ref{th2}}, and the universal algebra $\mathrm U(A,L)$ thus
acquires a Hopf algebra structure.

The case mentioned earlier where $A$ is the coordinate ring of an
algebraic group $H$ and where $\mathfrak g$ is the Lie algebra of
an algebraic group $G$ acting on $H$ by group automorphisms
is an example for this situation. More generally, $A$ could be
a general commutative Hopf algebra and   $\mathfrak g$ the Lie algebra of
an algebraic group $G$ acting on $A$ by Hopf algebra automorphisms.

Likewise, let $V$ be a rational representation of an algebraic
group $G$ over the ground field $\mathbf  k$. Addition in $V$
induces a Hopf algebra diagonal on the affine coordinate ring
$\mathbf  k[V]$ of $V$ that is compatible with the $G$-actions.
Consequently the familiar algebra of differential operators associated
with the representation then acquires a Hopf algebra structure over the
ground field.
\end{Example}

Thus, in the case at hand, we obtain a genuine Hopf algebra
structure on the universal algebra of the Lie-Rinehart algebra
under discussion.
It would be interesting to construct examples which are more
general than the above crossed products. Such examples arise,
perhaps, in differential Galois theory.

 In the papers \cite{kapraone} and
\cite{moermrcu}, for a general Lie-Rinehart algebra $(A,L)$ 
defined over a field of characteristic zero, it is
shown that
the classical construction of the diagonal map for the
universal enveloping algebra of an ordinary Lie algebra extends to
 an $A$-linear morphism from $\mathrm U(A,L)$ to $\mathrm
U(A,L) \otimes_A \mathrm U(A,L)$ which turns $\mathrm U(A,L)$ into
an $A$-coalgebra. However this does not yield a genuine bialgebra
structure unless $L$ is an ordinary $A$-Lie algebra so that
$\mathrm U(A,L)$ is then the ordinary universal algebra associated
with the $A$-Lie algebra $L$, and only a suitably defined
subalgebra $\mathrm U(A,L) \overline{\otimes}_A \mathrm U(A,L)$ of
$\mathrm U(A,L) \otimes_A \mathrm U(A,L)$ acquires an algebra
structure. The precise structure which  $\mathrm U(A,L)$ carries is that of an 
$A$-bialgebra, cf. \cite{kapraone} and the references there.
For intelligibility we recall that, given the commutative algebra $A$
 over the field
$\mathbf k$, a left $A$-bialgebra $H$ over 
$\mathbf k$ consists of an algebra $H$ (not necessarily commutative)
containing $A$, a morphism $\varepsilon\colon H\to A$ of left $A$-modules
which satisfies the identity
\[
\varepsilon(uv)= \varepsilon (u\cdot \varepsilon(v))
\]
together with a morphism $\Delta\colon H \to H \overline{\otimes}_A H$
of algebras into a suitably defined submodule
$H \overline{\otimes}_AH $
of the ordinary left $A$-module tensor product $H \overline{\otimes}_A H$,
the submodule $H \overline{\otimes}_AH $ being an algebra
under ordinary tensor product multiplication.

\section{Concluding remarks related with Lie-Rinehart bialgebras}

Let $L$ and $D$ be $(R,A)$-Lie algebras which, as $A$-modules, are
finitely generated and projective, in such a way that, as an
$A$-module, $D$ is isomorphic to $L^* = \mathrm{Hom}_A(L,A)$. We
say that $L$ and $D$ {\em are in duality\/}. We write $d$ for the
differential on $\mathrm{Alt}_A(L,A)\cong \Lambda_AD$ coming from
the Lie-Rinehart structure on $L$ and $d_*$ for the differential
on $\mathrm{Alt}_A(D,A)\cong \Lambda_AL$ coming from the
Lie-Rinehart structure on $D$. Likewise we denote the Gerstenhaber
bracket on $\Lambda_A L$ coming from the Lie-Rinehart structure on
$L$ by $[\, \cdot \, , \, \cdot \,]$ and that on $\Lambda_A D$
coming from the Lie-Rinehart structure on $D$ by $[\, \cdot \, ,
\, \cdot \,]_*$. The triple $(A,L,D)$ is said to constitute a {\em
Lie-Rinehart\/} bialgebra if the differential $d_*$ on
$\mathrm{Alt}_A(D,A) \cong \Lambda_A L$ and the Gerstenhaber
bracket $[\, \cdot \, , \, \cdot \,]$ on $\Lambda_A L$ are related
by
\[
d_*[x,y] = [d_*x,y] + [x,d_*y], \quad x,y \in L,
\]
or equivalently, if the differential $d$ on $\mathrm{Alt}_A(L,A)
\cong \Lambda_A D$ behaves as a derivation for the Gerstenhaber
bracket $[\, \cdot \, , \, \cdot \,]_*$ in all degrees, that is to
say
\[
d[x,y]_* = [dx,y]_* -(-1)^{|x|} [x,dy]_*, \quad x,y \in
\Lambda_AD.
\]
See \cite{banach} for details.

Let $(\Delta,\varepsilon)$ be
a bialgebra structure on $A$,  suppose that 
$(A,\Delta,L,\omega^{\otimes})$ is  a {\em bi-Lie-Rinehart
algebra\/}, and let
\begin{equation}
\Delta\colon \mathrm U(A,L) \longrightarrow \mathrm U(A,L)\otimes
\mathrm U(A,L) \label{diag3}
\end{equation}
and 
\begin{equation}
\varepsilon\colon \mathrm U(A,L) \longrightarrow R \label{coun3}
\end{equation}
be the resulting coalgebra structure on
the universal algebra $\mathrm U(A,L)$ 
given by Theorem \ref{th1} above
turning the latter into a bialgebra.
The diagonal map of $A$ induces a morphism
\begin{equation}
\mathrm{Hom}_A(L,A) \longrightarrow \mathrm{Hom}_A(L,A\otimes A)
\label{diag4}
\end{equation}
of $A$-modules.

Suppose that, as an $A$-module, $L$ is finitely 
generated and projective. Then
the canonical morphism
\begin{equation}
\mathrm{Hom}_A(L,A)\otimes  \mathrm{Hom}_A(L,A)  
\longrightarrow \mathrm{Hom}_{A\otimes A}(L\otimes L,A\otimes A)
\label{can5}
\end{equation}
is an isomorphism of $(A \otimes A)$-modules, and so is the canonical 
morphism
\begin{equation}
 \mathrm{Hom}_A(L,A\otimes A) \cong
 \mathrm{Hom}_{A\otimes A}(L^{\otimes},A\otimes A)  .
\label{can4}
\end{equation}
Let
\begin{equation}
[\,\cdot\, ,  \, \cdot \,]\colon \mathrm{Hom}_A(L,A) \otimes
 \mathrm{Hom}_A(L,A) \longrightarrow \mathrm{Hom}_A(L,A) 
\label{bra}
\end{equation}
be a Lie bracket on $\mathrm{Hom}_A(L,A)$  turning $(A,L,  \mathrm{Hom}_A(L,A) )$
into a Lie-Rinehart bialgebra.
The composite of \eqref{bra} with the inverse of the isomorphism \eqref{can5}
and with the morphism \eqref{diag4} 
yields the morphism
\begin{equation}
 \mathrm{Hom}_{A\otimes A}(L\otimes L,A\otimes A)
 \longrightarrow  \mathrm{Hom}_{A}(L,A\otimes A)  
\label{can6}
\end{equation}
which, in view of the isomorphism \eqref{can4}, takes the form
\begin{equation}
 \mathrm{Hom}_{A\otimes A}(L\otimes L,A\otimes A)
 \longrightarrow  \mathrm{Hom}_{A\otimes A}(L^{\otimes},A\otimes A).  
\label{can7}
\end{equation}
The $(A\otimes A)$-dual thereof dual---here we use the hypothesis that,
as an $A$-module, 
$L$ is finitely generated and projective---is a morphism of the kind
\begin{equation*}
L^{\otimes}\longrightarrow L\otimes L .
\end{equation*}
The composite thereof with the morphism
$L\to L^{\otimes}$
which is a constituent of the morphism
\begin{equation}
(\Delta,\Delta \otimes \mathrm{Id}) \colon (A,L) \longrightarrow
(A\otimes A,L^{\otimes}) \label{ob44}
\end{equation}
of Lie-Rinehart algebras given above as \eqref{ob4}
yields a morphism of the kind
\begin{equation}
L\longrightarrow L\otimes L .
\label{pert}
\end{equation}
We conjecture that perturbing the comultiplication
\eqref{diag3} by means of the morphism \eqref{pert}
yields a new 
comultiplication which, together with the other data, again combines to a
bialgebra structure on 
$\mathrm U(A,L)$.
Since we do not have a particular example yet,
we will pursue these issues elsewhere.

\end{document}